\documentclass[12pt]{article}
\usepackage{amssymb,amsmath}

 \topmargin= - 8 true mm
\newcommand{\pr}{\mathbb{P}}                     
\newcommand{\exn}{\mathbb{E}\,}                  
\newcommand{\R}{\mathbb{R}}                      
\newcommand{\deq}{\stackrel{d}{=}}               

\newtheorem{theorem}{Theorem}

\begin{document}

\title{On limiting cluster size distributions for processes of exceedances for stationary sequences}

\author{K. Borovkov\footnote{Department of Mathematics and Statistics,
The University of Melbourne, Parkville VIC 3010, Australia. E-mail:
K.Borovkov@ms.unimelb.edu.au. Research supported by ARC Grant DP0880693.} \ and S.
Novak\footnote{Middlesex University, MUBS, The Burroughs, London NW44BT, UK. E-mail:
S.Novak@mdx.ac.uk. Supported by the EPSRC grant EP/H006613.}}

\date{}

\maketitle

\begin{abstract}
It is well known that, under broad assumptions, the time-scaled point process of
exceedances of a high level by a stationary sequence converges to a compound Poisson
process as the level grows. The purpose of this note is to demonstrate that, for any
given distribution $G$ on $\mathbb{N}$, there exists a stationary sequence for which
the compounding law of this limiting process of exceedances will coincide with~$G$.


\end{abstract}

\bigskip

Let $\{X_k\}_{k\ge 0}$ be a stationary discrete time real-valued process. For a suitably
chosen increasing real sequence $\{u_n\}$, consider the time-scaled point process of
exceedances
\begin{equation}
N_n (A):= \#\{k/n\in A: \, X_k > u_n\}, \qquad A\in \mathcal{B}(\mathbb{R}_+).
 \label{PP}
\end{equation}
As is well known (see e.g. Corollary~3.3 and Theorem~4.1 in~\cite{HsHuLe} and
also~\cite{Ro}), under broad assumptions, if the process $N_n$ has a limit as
$n\to\infty$, the latter must be a compound Poisson process, with some compounding
law~$G=\{g_k\}_{k\ge 1}$ on $\mathbb{N}$. An example of an important class of processes
for which $G$ is non-trivial can be found e.g.\ in~\cite{BoLa} (Scenario~4.3): in the
case of asymptotically homogeneous Markov chains, $G$ will be geometric.

The main objective of this short note is to complete  the characterization of the class
of limiting distributions for~\eqref{PP} by giving an affirmative answer to the
following natural question: {\em Given an arbitrary distribution $G$ on $\mathbb{N}$,
does there exist a stationary process $\{X_\bullet\}$ for which $G$ will be the cluster
size distribution for the limiting point process of exceedances?} This is achieved by
constructing a two-dimensional stationary Markov chain (essentially, a regenerative
process), then taking $\{X_\bullet\}$ to be the component process of the chain and
applying to it results from~\cite{Ro}. While in the case of a finite mean
\[
\mu :=\exn \zeta, \qquad \zeta\sim G,
\]
this task is next to trivial, the case $\mu=\infty$ is more interesting and is, in our
opinion,  worth presenting.

\begin{theorem}
For any distribution $G$ on $\mathbb{N}$ there exists a stationary process
$\{X_\bullet\}$ for which the time-scaled point process of exceedances~\eqref{PP}
converges, for a suitably chosen $\{u_n\uparrow\}$, to a compound Poisson process with
compounding law~$G$.
\end{theorem}

{\em Proof.} Without loss of generality, we can assume that
\begin{equation}
\mbox{g.c.d.}\{k\ge 1: \, g_k=1\}=1.
 \label{1}
\end{equation}
In the case when $\mu < \infty$, the construction is straightforward: it will be a
regenerative process staying at randomly chosen levels during regeneration cycles of
random lengths distributed according to~$G$, the heights of the levels and the lengths
of the cycles forming independent i.i.d.\ sequences. Then exceedances of a high level
will automatically be clustered, with cluster size distributed almost as~$\zeta$ since
it is quite unlikely to have two cycles with large heights one after another.

More formally, let $\{Y_k\}_{k\ge 1}$ be an i.i.d.\ sequence of random variables that
are, say, exponentially distributed with mean one, and let $\{\tau_k\}_{k\ge 1}$ be a
sequence of independent random variables which is independent of~$\{Y_\bullet\}$ and
such that
\begin{equation}
\tau_n \deq \zeta, \quad n\ge 2; \qquad \pr (\tau_1 = j)=
 \frac1\mu\, \pr (\zeta \ge j)\equiv  \frac1\mu\sum_{m\ge j} g_m,
\quad j\ge 1.
 \label{TA}
\end{equation}
Putting $S_k:= \sum_{j=1}^k \tau_j,$ we see that
\[
\eta(n) :=\min\{k\ge 1:\,   S_k> n\}
\]
is a (delayed) stationary renewal process, with a linear renewal function:
\begin{equation}
\sum_{j=1}^\infty \pr (S_j=k) = \frac1{\exn \tau_2} =\frac1\mu, \quad k\ge 1
 \label{RF}
\end{equation}
(see e.g.\ Section~9.2 in~\cite{Bo}). Therefore the process
\begin{equation}
X_k:= Y_{\eta(k)}, \quad k\ge 0,
 \label{SPP}
\end{equation}
will also be stationary: $\{X_k\}_{k\ge 0}\deq \{X_{l+k}\}_{k\ge 0}$ for any $l\ge 1$,
which immediately follows from the independence of the sequences $\{\tau_\bullet\}$ and
$\{Y_\bullet\}$ and the well-known fact that, for any $l\ge 1$, one has $S_{\eta (l)}-l
\deq \tau_1$.

The process $\{X_\bullet\}$ is clearly regenerative in the sense of~\cite{Ro}, with
i.i.d.\ cycles $C_k:=\{X_{S_{k-1}+j}:\, 0\le j< S_{k} - S_{k-1}\}$, $k>1$. It is
obvious that, for the number of exceedances of $\{X_\bullet\}$ of the level $u_n$
during the cycle $C_k$ defined as
\[
\xi_k := \# \{j :\, S_{k-1} \le j < S_{k}, \ X_j > u_n\},
\]
one has
\begin{equation}
\xi_k = \left\{
\begin{array}{ll}
   \tau_k, & \mbox{if}\  Y_k > u_n,\\
 0, & \mbox{otherwise},
 \end{array}
 \right.
 \label{GS-}
\end{equation}
and so
\begin{equation}
 \pr (\xi_k = j \, |\,\xi_k >0 ) =g_j, \quad k\ge 2.
 \label{GS}
\end{equation}
Now it follows immediately from Theorem~3.3 in~\cite{Ro} that $G$ will be the asymptotic
distribution of the exceedances cluster size.

Now turn to the case $\mu=\infty.$ The simple construction presented above won't work
in this case as the regeneration cycles will have infinite mean lengths, but it can be
modified by making the components of the random vectors $(\tau_k, Y_k),$ $k>1$,
dependent of each other in such a way that (i)~the conditional distribution of $\tau_k$
given $Y_k=y$ converges to~$G$ as $y\to\infty$ (as we are interested in exceedances of
high levels after all, we need to control the cycle length law only when there is an
exceedance inside the cycle) and (ii)~$\exn \tau_k<\infty$. Observe that a different
regenerative process with the cycle length distribution depending on the level height
was used in~\cite{Sm} to give a counterexample concerning the interpretation of the
extremal index, see Remark~2 below.

To construct our version of the modified regenerative process, consider   an i.i.d.\
sequence $\{(\zeta_k, Y_k)\}_{k\ge 2}$ of random vectors with independent components,
$\zeta_k\deq\zeta$, $Y_k$ being exponential random variables with unit mean, and let
$\tau_k := \zeta_k \wedge \lceil Y_k\rceil ,$ $k\ge 2 ,$ where, as usual, $\lceil
x\rceil: = \min\{k\in \mathbb{N} :\, k\ge x\}$.

Thus defined i.i.d.\ random variables $\tau_k$ will have the distribution
\[
p_j := \pr (\tau_2 = j) = \int_{0}^\infty  g_{j} (v) e^{-v} dv, \quad j\ge 1,
\]
where
\[
g_{j} (v) := \pr ( \zeta \wedge  \lceil v\rceil =j) = \left\{
 \begin{array}{ll}
 g_j, & j< \lceil v\rceil,\\
 \overline{g}_{\lceil v\rceil}, & j=\lceil v\rceil,\\
 0, & j > \lceil v\rceil,
 \end{array}
\right. \qquad  \overline{g}_m:=  \sum_{i\ge m} g_i,
\]
with the   mean
\begin{equation}
\nu := \exn \tau_2 =\int_0^\infty \exn (\zeta \wedge \lceil v\rceil )  e^{-v} dv
 \le \int_0^\infty   \lceil v\rceil    e^{-v} dv
<\infty.
 \label{Tau2Mean}
\end{equation}
Observe that, for $v>0$
\begin{equation}
\pr (Y_2 \in dv  \, |\, \tau_2 =  j)
   =   \frac{\pr (\tau_2 =  j \, |\, Y_2 = v  )}{\pr (\tau_2 = j)}
   \, \pr( Y_2 \in dv  )
   = \frac{g_{j}(v)}{p_j} \, e^{-v} dv.
 \label{TauCon}
\end{equation}

Now we have to define $(\tau_1, Y_1)$ (assumed to be independent of $\{(\tau_k,
Y_k)\}_{k\ge 2}$) in such a way that the process~\eqref{SPP} will again be stationary.
It is obvious that, similarly to~\eqref{TA}, $\tau_1$ should follow the distribution
\begin{equation}
\pr (\tau_1 = j)=
 \frac1\nu\, \pr (\tau_2 \ge j)\equiv  \frac1\nu\sum_{m\ge j} p_m, \quad j\ge 1,
 \label{TAt}
\end{equation}
so we only need to specify how $Y_1$ depends on $\tau_1$. Again, it's quite clear
that the dependence should be the same as one has in the limit (as $k\to\infty$)
between $Y_{\eta(k)}$ and the overshot~$S_{\eta(k)} -k$. As the length $\tau_{\eta(k)}$
of the renewal interval `covering' the point~$k$ is, loosely speaking, greater than
that of a `typical' $\tau,$ our construction implies that $Y_{\eta(k)}$ should also be
greater than a `typical'~$Y$ (and, in particular, cannot have the same distribution
as~$Y_2$, cf.~\eqref{Y1}). The above informal argument leads to the following
construction.

Denote by
\[
\gamma (k) := k - S_{\eta(k)-1}, \qquad \chi (k) := S_{\eta(k)} -k
\]
the defect and excess of the level $k$ in the random walk $\{S_{\bullet}\}$,
respectively. As is well known (recall~\eqref{1} and see e.g.\ Section~9.3
in~\cite{Bo}),
\begin{equation}
\lim_{k\to\infty} \pr (\gamma (k) =i, \, \chi (k)=j) = \frac{p_{i+j}}{\nu}, \qquad i\ge
0, \ \ j\ge 1.
 \label{GX}
\end{equation}
Now consider a random vector $(\gamma, \chi, V)$ assuming that its first two components
are integer-valued and such that $\pr (\gamma  =i, \, \chi  =j) $ is given by the RHS
of~\eqref{GX} (note that the distribution of~$\chi$ will coincide with that of $\tau_1$
from~\eqref{TAt}), whereas
\[
\pr (V\in dv \, |\, \gamma  =i, \, \chi  =j) = \pr (Y_2 \in dv  \, |\, \tau_2 = i+j),
\qquad i\ge 0, \ \ j\ge 1.
\]
Finally, we set
\[
(\tau_1, Y_1)\deq (\chi, V)
\]
and again consider the process $\{X_\bullet\}$ defined by~\eqref{SPP}.

It is obvious from the construction that, to prove that $\{X_\bullet\}$ is stationary,
it suffices to show that, for any $l\ge 1$, $j\ge 1$, $m\ge 1$ and $v>0$, one has
\begin{equation}
\pr (Y_{\eta(l)} \in dv, \, S_{\eta(l)}-l =j ) = \pr (Y_{1} \in dv, \, \tau_1 =j ).
 \label{LR}
\end{equation}
To do that, we first observe that
\begin{align}
\pr (Y_{1} \in dv, \, \tau_1 =j )
  & = \pr (V\in dv,   \,   \chi  =j) \notag \\
  & = \sum_{i\ge 0} \pr (V\in dv  \, |\, \gamma  =i, \, \chi  =j)\, \pr( \gamma =i, \, \chi
  =j) \notag \\
  & =  \sum_{i\ge 0} \pr (Y_2 \in dv  \, |\, \tau_2 = i+j) \frac{p_{i+j}}{\nu}\notag \\
  & = \frac{e^{-v} dv}{\nu} \sum_{i\ge 0} g_{ i+j} (v)
  \label{LRa}
\end{align}
from~\eqref{TauCon}. Now, for the LHS of~\eqref{LR}, we have
\begin{align}
\pr (Y_{\eta (l)} \in dv,   \, S_{\eta(l)}-l =j )
   & = \pr (Y_1 \in dv, \, \tau_1=l+j)  \notag\\
  &  \quad +
     \sum_{i=0}^{l-1} \sum_{r=1}^\infty \pr (Y_{r+1} \in dv, \, S_r = l-i, \, S_{r+1} =
     l+j)\notag \\
     & = \frac{e^{-v} dv}{\nu} \sum_{i\ge 0} g_{l+i+j} (v)\notag \\
  &  \quad +
     \sum_{i=0}^{l-1} \sum_{r=1}^\infty
     \pr (S_r = l-i)\, \pr (Y_{r+1} \in dv, \, \tau_{r+1} =
     i+j),
     \label{XA}
\end{align}
where we used \eqref{LRa} to evaluate the first   term on the RHS of the first line.
The inner sum in   the last line of~\eqref{XA},  using the first relation in~\eqref{RF}
and~\eqref{TauCon}, is equal to
\begin{multline*}
\pr (Y_{2} \in dv  ,\, \tau_{2} = i+j) \sum_{r=1}^\infty  \pr (S_r = l-i) \\
= \frac{\pr (Y_{2} \in dv  \, |\,  \tau_{2} = i+j)}{\exn\tau_2} \, \pr (\tau_2= i+j)
 = \frac{e^{-v}dv}{\nu} \, g_{i+j} (v).
\end{multline*}
Substituting these expressions into the RHS of~\eqref{XA} yields
\[
\pr (Y_{\eta(l)} \in dv,   \, S_{\eta(l)}-l =j )
 = \frac{e^{-v} dv}{\nu} \left[
  \sum_{i\ge 0} g_{l+i+j} (v) + \sum_{i=0}^{l-1} g_{i+j} (v)
  \right] = \frac{e^{-v} dv}{\nu} \sum_{i\ge 0} g_{ i+j} (v),
\]
which coincides with the RHS of \eqref{LRa} thus proving the desired stationarity
of~$\{X_\bullet\}$.

It remains to observe that our process  $\{X_\bullet\}$ is again regenerative, with the
finite mean cycle length~$\nu$, that~\eqref{GS-} still holds for it, and that, instead
of~\eqref{GS}, we   now have
\[
 \sup_{j\ge 1}|\pr (\xi_k = j \, |\,\xi_k >0 ) - g_j| \le \overline{g}_{\lceil
 u_n\rceil}.
\]
Clearly, $\overline{g}_{\lceil u_n\rceil}\to 0$ as $n\to\infty$, and so condition (3.9)
of Theorem~3.3 in~\cite{Ro} will be satisfied. Therefore the claim of the theorem will
hold in this case as well, showing that $G=\{g_j\}$ will emerge as the asymptotic
distribution for the size of clusters of exceedances. Theorem~1 is proved.

\hfill$ \square$

\medskip\noindent{\bf Remark~1} Observe that it follows from \eqref{LRa} and
\eqref{Tau2Mean} that
\begin{equation}
 \label{Y1}
\pr (Y_1 \in dv) = \frac{e^{-v} dv}{\nu} \sum_{j\ge 1}\sum_{i\ge 0} g_{ i+j} (v)
 = \frac{\pr (Y_2 \in dv)}{\exn \tau_2} \sum_{j\ge 1} \pr (\tau_2\ge j\,|\, Y_2=v).
\end{equation}
As, due to stationarity   $X_k\deq X_1 = Y_1$ for any $k\ge 1$, the above implies that
the distribution of $X_k$ in our stationary process is given by
\[
\pr (X_k\in dv) \equiv \pr (Y_1 \in dv) = \frac{\exn (\tau_2 \,|\, Y_2 =v)}{\exn \tau_2}
\,
 \pr (Y_2 \in dv), \quad v>0.
\]

\medskip\noindent{\bf Remark~2} Recall that, under broad assumptions, the
limiting distributional type for   the maxima in a stationary sequence coincides with
that for the maxima of i.i.d.\ random variables with the same marginal distribution,
and that the changes brought by dependence can often be characterized by the so-called
``extremal index" of the stationary sequence (see e.g.\ Section~2 in~\cite{LeRo} and
further references therein, and also~\cite{No}). For a  sequence
$\{x_\bullet\}\subset\R$, set $M_n(x_\bullet):=\max_{1\le k\le n} x_k$, and let $\{\hat
X_\bullet\}$ be an i.i.d.\ sequence with $\hat X_k \deq X_k.$ As it is well known (see
e.g.\ Lemma~1.2.2 in~\cite{LeRo}), for a sequence of constants $\{u_\bullet \}$ and
$\lambda \in [0,\infty]$, one has $n \pr (X_0 > u_n) \to \lambda$ as $n\to\infty$ iff
\begin{equation}
 \label{Mo}
 \pr (M_n (\hat X_\bullet) \le u_n) \to e^{-\lambda}.
\end{equation}
It turns out that, in many cases (originally it was noted for strongly mixing sequences
in~\cite{Lo}),  if~\eqref{Mo} holds then one also has
\[
\pr (M_n ( X_\bullet) \le u_n) \to e^{-\theta\lambda}
\]
for a some fixed value $\theta\in [0,1]$ which is referred to as the {\em extremal
index\/} of the sequence $\{X_\bullet\}$.

Alternatively, the extremal index can be characterized by the fact that $1/\theta$  is
the limiting mean cluster size in the sequence of point processes~\eqref{PP} of
exceedances over high levels, as it was shown under broad assumptions in~\cite{HsHuLe};
a counterexample showing that this interpretation of~$\theta$ is not necessarily
correct was given in~\cite{Sm}.  For the stationary sequences that we constructed in
the proof of Theorem~1, the extremal index is equal to $\theta =1/\exn \zeta$ in both
cases ($\theta =0$ when $\exn \zeta=\infty$), which can easily be verified by a direct
calculation.

\end{document}